\documentclass[11pt,lefteqn]{article}

\usepackage{amsmath,amsfonts,amsthm,amssymb,graphics,epsfig,color}
\newcommand{\R}{\mathbb R}

\newcommand{\D}{\mathbb D}

\renewcommand{\Re}{\mathop {\rm Re}\nolimits}
\renewcommand{\Im}{\mathop {\rm Im}\nolimits}

\begin{document}
\date {\today}
\title{Numerical bounds on the Crouzeix ratio for a class of matrices}

\author{Michel Crouzeix\footnote{Univ.\,Rennes, CNRS, IRMAR\,-\,UMR\,6625, F-35000 Rennes, France.
email: michel.crouzeix@univ-rennes.fr},  Anne Greenbaum\footnote{University of Washington,
Applied Math Dept., Box 353925, Seattle, WA 98195.  email:  greenbau@uw.edu},
{Kenan Li}
}

\maketitle

\begin{abstract}  We provide numerical bounds on the Crouzeix ratio
for KMS matrices $A$ which have a line segment on the boundary 
of the numerical range. The Crouzeix ratio is the supremum 
over all polynomials $p$ of the spectral norm of $p(A)$ divided
by the maximum absolute value of $p$ on the numerical range of $A$.
Our bounds confirm the conjecture that this ratio
is less than or equal to $2$. We also give a precise description of these numerical ranges.
\end{abstract}

\paragraph{2000 Mathematical subject classifications\,:} 15A60; 15A45; 47A25 ; 47A30

\noindent{\bf Keywords\,:}{ numerical range, spectral set}

\section{Introduction}
We consider $n$ by $n$ matrices $A_n$ with $1$'s in the strict
upper triangle and $0$'s elsewhere.  For $n=3,4,5,6$, we have numerically
determined upper and lower bounds on the value
\[
\psi (A_n) := \sup \{ \| p( A_n ) \| : \mbox{$p$ a polynomial with $|p| \leq 1$
in $W( A_n )$} \} ,
\]
where $W( A _n )$ denotes the numerical range of $A_n$ and $\| \cdot \|$ is the
spectral norm in ${\mathbf C}^{n,n}$. (We refer to this quantity 
as the {\em Crouzeix ratio} although sometimes this term denotes the reciprocal,
$\max_{z \in W( A_n )} |p(z)| / \| p( A_n ) \|$, for a given polynomial $p$ \cite{Overton2022}.)
Crouzeix's conjecture is that for
all square matrices $A$, $\psi (A) \leq 2$.  A way to determine $\psi (A)$ is to
introduce a Riemann mapping $g$ from the interior of $W(A)$ onto the
open unit disk $\mathbb{D}$ and to consider the matrix $M := g(A)$.
In this case, we can write
\[
\psi (A) = \psi_{\mathbb{D}} (M) := \max \{ \| f(M) \| : \mbox{$f$ holomorphic
in ${\mathbb{D}}$ with $|f| \leq 1$ in ${\mathbb{D}}$} \} .
\]
We know that the maximum is realized by a Blaschke product of order $n-1$
and each choice of such a Blaschke product $b$ provides a lower bound:
$\psi (A) = \psi_{\mathbb{D}} (M) \geq \| b(M) \|$.  From the von Neumann
inequality, we easily deduce 
\[
\psi_{\mathbb{D}} (M) \leq \psi_{cb,{\mathbb{D}}} (M) := \min \{ \mbox{cond} (H) : 
H \in \mathbb{C}^{n,n} ,~ \| H^{-1} M H \| \leq 1 \} .
\]
Thus, the exhibition of a matrix $H$ satisfying $\| H^{-1} M H \| \leq 1$
leads to an upper bound $\psi (A) \leq \mbox{cond} (H) := \| H \| \cdot 
\| H^{-1} \|$.  This approach has been used to prove $\psi (A) \leq 2$
for some classes of matrices \cite{Crouzeix2004,Choi2013,ChoiGreen2015,
GKL2018}, but this assumes that one knows the matrix $M$ with sufficient
accuracy.  This is the case when $W(A)$ is an ellipse since there is an
analytic formula for $g$ and also for \cite{Choi2013,ChoiGreen2015} where
$g(A) = cA$ for a certain constant $c$.  In general, however, there is no
simple expression for the boundary of $W(A)$ and for the Riemann mapping $g$.
There are numerical methods for computing $g$ and thus $M$ with high precision,
but to {\em guarantee} the accuracy would require a complete analysis
of all discretization and rounding errors.  

In section 2, we consider the matrix 
\[
A_3 = \left[ \begin{array}{ccc} 0 & 1 & 1 \\ 0 & 0 & 1 \\ 0 & 0 & 0 \end{array}
\right] .
\]
Our numerical and analytical work suggests that $1.9956978 < \psi ( A_3 ) <
1.9956979$, and we are confident in this range, although it does not
provide a proof that $\psi ( A_3 ) \leq 2$, since it relies on numerical
computation of $g(A)$.

Another approach is to identify a rational function $f_1$ such that the
image $\Omega := f_1 ( \mathbb{D} )$ of the unit disk is a subset of
(and close to) $W( A_3 )$.  Let $g_1$ be the inverse of $f_1$.  Then with
$M_1 = g_1 ( A_3 )$, we can write
\[
\psi_{\Omega} ( A_3 ) := \sup \{ \| h( A_3 ) \| : |h| \leq 1 \mbox{ in } \Omega \}
= \psi_{\mathbb{D}} ( M_1 ) := \sup \{ \| f( M_1 ) \| : |f| \leq 1 
\mbox{ in } \mathbb{D} \} .
\]
Since $\Omega \subset W( A_3 )$, we clearly have $\psi ( A_3 ) \leq 
\psi_{\Omega} ( A_3 ) \leq \psi_{cb,\mathbb{D}} ( M_1 )$.  We are now 
able to compute $M_1$ analytically and exhibit a matrix $H_1$ such
that $\| H_1^{-1} M_1 H_1 \| = 1$ and $\mbox{cond} ( H_1 ) \approx 1.9999514$.
However, we must verify numerically that $\Omega \subset W( A_3 )$.

In section 3, we consider more generally the matrices 
$A_n$ for $n>3$. They belong to the class of KMS matrices \cite{KMS1953} and they are the matrices in this class for which the boundary
of the numerical range
contains a line segment \cite{GW2013}. 
We derive a simple description of their
numerical ranges and determine numerically the following bounds:
\[
1.993800 \leq \psi ( A_4 ) \leq 1.993801 ,~~
1.992921 \leq \psi ( A_5 ) \leq 1.992922 ,~~
1.992444 \leq \psi ( A_6 ) \leq 1.992445 .
\]

In section 4, we explain the numerical method used to compute the conformal
mapping and we provide the Matlab code used for the computation of $M = g(A)$.

\section{Numerical estimates for the matrix $A_3$}

Here we consider the matrix
$A_3=\begin{pmatrix}0 &1 &1\\0&0&1\\0&0&0\end{pmatrix}$.
 We will see in the next section that the boundary of its numerical range
 is the union of a part of an algebraic curve
$\Big\{ \frac{2\,e^{i\theta }+e^{2i\theta }}{3}\,: -\frac{2\pi }{3}\leq \theta \leq \frac{2\pi }{3}\Big\}$
and of the vertical straight line $[-\frac12-i\frac{\sqrt3}{6},-\frac12+i\frac{\sqrt3}{6}]$.
The algebraic curve is a cardioid\,; its Cartesian equation is
$
27(x^2{+}y^2)^2-18(x^2{+}y^2)-8x-1=0.
$
\begin{figure}[h]
\centerline{\psfig{figure=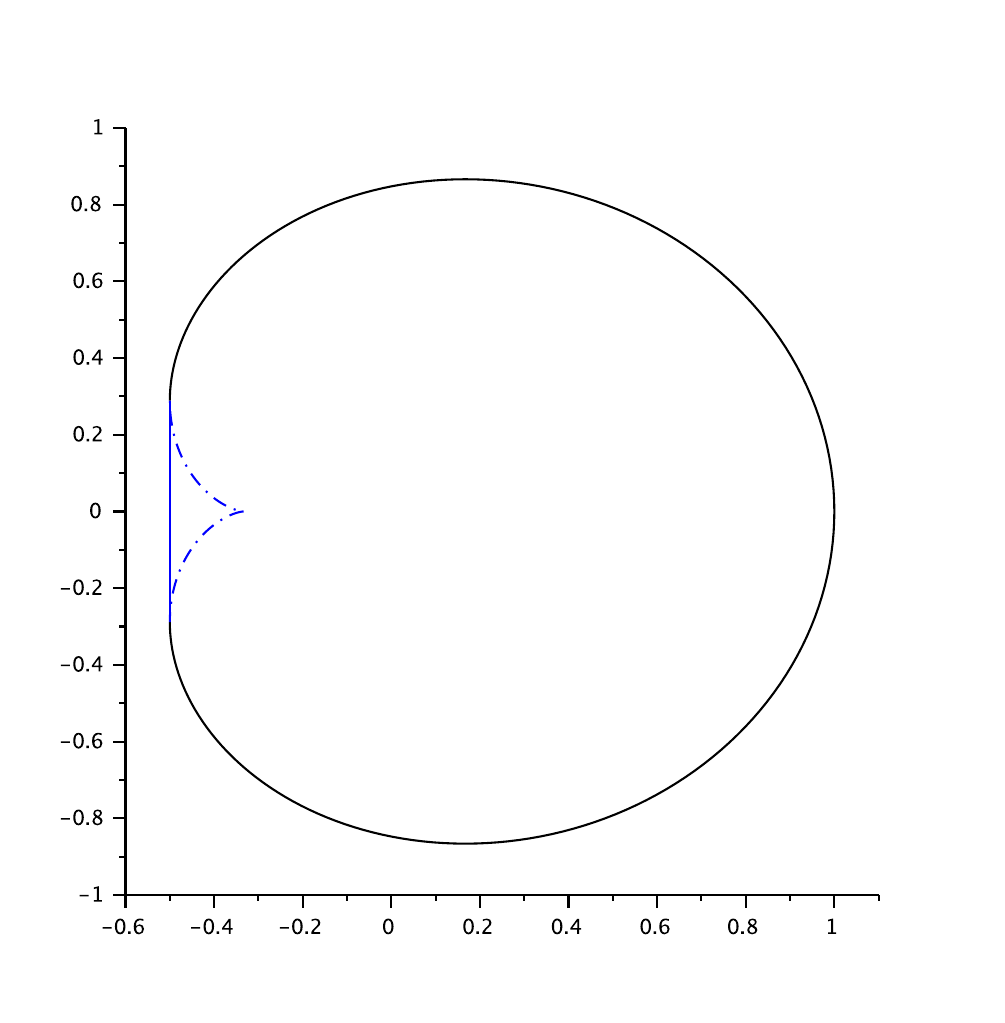 ,height=7cm,width=7cm}}
\caption{{ The boundary of the numerical range in black, the remaining part of the algebraic curve in dashed blue}}\label{fig1}
\end{figure}
\bigskip

We denote by $g$ the Riemann mapping from $W(A)$ onto the unit disk $\D$ that satisfies $g(0)=0$ and $g'(0)>0$. Then 
\begin{align*}
M=g(A)=\begin{pmatrix}
0 &a &b\\0&0&a\\0&0&0
\end{pmatrix},\quad \textrm{with }
a=g'(0),\ b=g'(0)+\tfrac12g''(0).
\end{align*}
From the numerical computations it appears that $a\simeq1.360374515$, $b\simeq0.710915425$ with an accuracy that we empirically estimate better than $10^{-8}$.
Using the Blaschke product $\displaystyle f(z)=\frac{z+0.5470208}{1+0.5470208z}\ \frac{z-0.1465739}{1-0.1465739z}$, we obtain
$\|f(M)\|=1.9956978$, which (numerically) shows that $\psi (A)\geq 1.9956978$.

We now choose the matrix
\begin{align*}
H=\begin{pmatrix}
a &b/2a &-b^2/8a^3\\0&1&-b/2a^2\\0&0&1/a
\end{pmatrix}, \quad\hbox{then}\quad
H^{-1}MH=\begin{pmatrix}
0&1&0\\0&0&1\\0&0&0
\end{pmatrix}.
\end{align*}
Therefore
$\psi_{cb,\D}(M)\leq \|H\| \|H^{-1}\|=$\,cond$(H) \simeq 1.995697855$.

With the numerical values obtained previously, we believe that we have the two sided estimates $1.9956978<\psi_M(a,b)<1.9956979$.
Therefore, it appears from the numerical simulation that $W(A)$ is a (complete) $1.9956979$-spectral set for $A$, that the complete bound\,\cite{paul} is the same as the ordinary bound, and that a function which realizes $\psi (A)$ is a Blaschke product of order 2 with 2 real roots.
But we have only an empiric estimate of the accuracy that we justify as follows. 
If we let $a(n)$ be the numerical value of $a$ computed with our program (described further) using $n$ points on the boundary, we have verified that $33\leq \frac{a(1447)-a(n)}{n^{-4}}\leq 65$ for values of $n$ between $23$ and $1205$. This suggests that our method is of order $n^{-4}$ and suggests that $|a(1205)-a|\leq 2 \cdot 10^{-11}$. Similarly, it appears that our computation of $b$ is of order $n^{-4}$, $130\leq \frac{b(n){-}b(1447)}{n^{-4}}\leq 350$ and that $|b(1205)-b|\leq 10^{-10}$.\bigskip

We turn now to the second attempt which is to consider the image $\Omega=f_1(\D)$ of the unit disk by the rational function 
$f_1(z)=(c_1 z+c_2 z^2+\cdots+c_7 z^7 )/(1+d_1 z+\cdots+d_7 z^7 )$, with the values 
\begin{align*}
c&=(0.734,0.49736,0.07268,-0.00521,0.00013,0.00061,-0.00251),\\
d&=(0.32564,-0.03291,0.01,-0.004,0.00084,-0.00242,0.00028).
\end{align*}
We let $g_1$ be the inverse function of $f_1$ and we set $M_1=g_1(A)$. We will see that if $\Omega$ is included in $W(A)$, then 
\[
\psi (A)\leq \psi_\Omega(A)\leq \psi _{cb,\D}(M_1)\leq {\rm cond}(H_1),\quad{\rm if \ \ }\|H_1^{-1}M_1H_1\|\leq 1.
\]
\begin{figure}[!ht] \centering
\begin{minipage}[t]{8cm} \centering
\includegraphics[width=8cm]{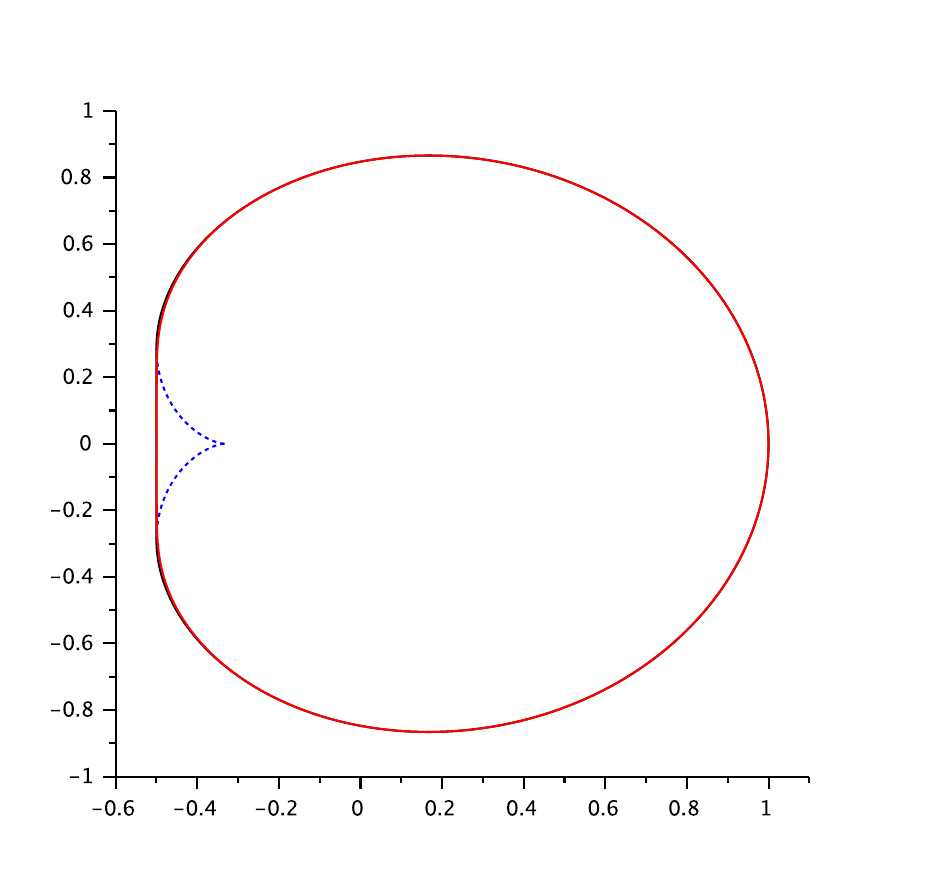}
\caption{The boundary of $\Omega$ in red, of $W(A)$ in black.}\label{fig1}
\end{minipage} \hskip1cm
\begin{minipage}[t]{6cm}
\centering \includegraphics[width=3cm]{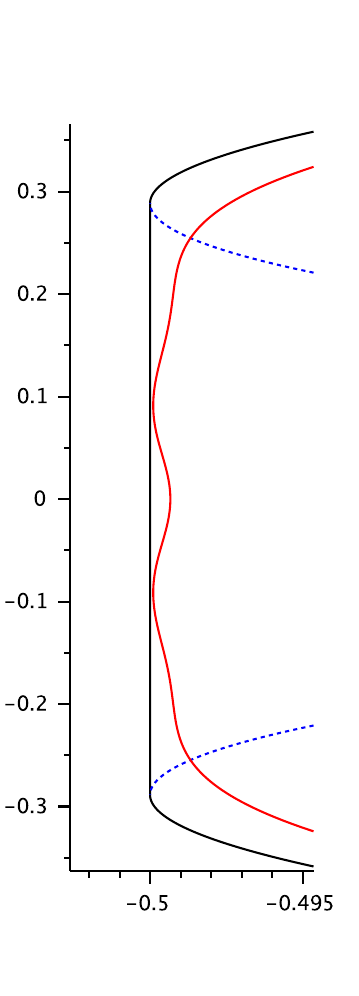} 
\caption{Zoom close to the straight line.}
\end{minipage} \end{figure}

 We are now able to compute $a_1=g_1'(0)$, $b_1=g_1'(0)+\frac12g_1''(0)$, and thus $M_1$ with an accuracy better than $10^{-14}$. We choose  
 \begin{align*}
H_1=\begin{pmatrix}
a_1 &b_1/2a_1 &-b_1^2/8a_1^3\\0&1&-b_1/2a_1^2\\0&0&1/a_1
\end{pmatrix}, \quad\hbox{then}\quad
H_1^{-1}M_1H_1=\begin{pmatrix}
0&1&0\\0&0&1\\0&0&0
\end{pmatrix}.
\end{align*}
This gives an estimate $\psi (A)\leq \psi_{cb,\D}(M_1)\leq 1.9996222$ with an accuracy better than $10^{-12}$ which ensures that $W(A)$ is a 2-spectral set for $A$. \bigskip

It remains to show that $W(A)$ contains $\Omega$.
For that, we first remark that the set $\{ z\,:p(z)<0\}$
with $p(z):=27|z|^4-18|z|^2-8\Re z-1<0$ is the interior of the cardioid and that the rectangle $\{z\,: - \tfrac12\leq \Re z\leq0$ and $|\Im z|\leq \sqrt{3}/6\}$ is contained in $W(A)$.
Taking into account the symmetry with respect to the real axis, in order to show that $\Omega$ is contained in $W(A)$, it suffices to show that the set
$\{z=f_1(e^{i\theta})\,: 0\leq \theta \leq \frac{3\pi }{4}\}$ is interior to the cardioid and that the set $\{z=f_1(e^{i\theta})\,: \frac{3\pi }{4}\leq  \theta \leq\pi\}$ is interior to the rectangle. \medskip

a) With $\theta _j=\frac{j\pi }{1000}$, $0\leq j\leq 750$, we have computed 
 $\displaystyle\max_{0\leq j\leq 750}p(f_1(e^{i\theta_j }))=-0.0008777\dots$ and 
$\displaystyle\max_{0\leq j\leq 749}\Big|\frac{p(f_1(e^{i\theta_{j+1 }}))-p(f_1(e^{i\theta_{j }}))}{\theta _{j+1}-\theta _j}\big|=0.0174\dots$ This gives us, for $0\leq \theta \leq 3\pi /4$, an estimate of $\max\big|\frac {d}{d\theta }p(f_1(e^{i\theta }))\big|\leq 0.018$ and thus $\max p(f_1(e^{i\theta}))\leq -0.0008777+0.018\pi /2000<-0.000849$. This shows that the set $\{z=f_1(e^{i\theta})\,: 0\leq \theta \leq \frac{3\pi }{4}\}$ is interior to the cardioid, thus interior to $W(A)$.
\begin{figure}[!ht] \centering
\begin{minipage}[t]{6cm} \centering
\includegraphics[width=8cm]{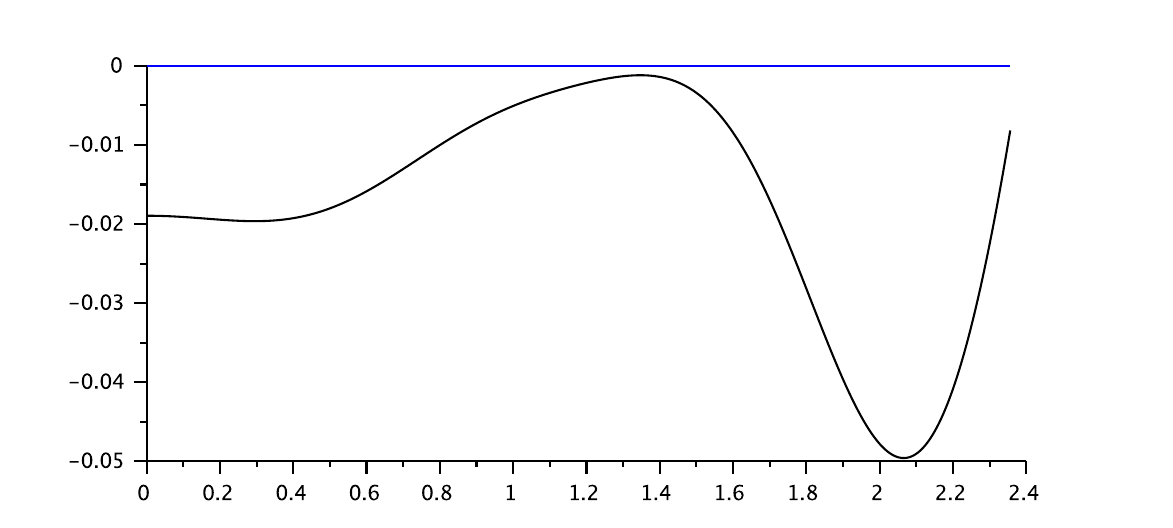}
\caption{Curve $p(f_1(e^{i\theta}))$, $0\leq \theta \leq\frac{3\pi }{4}$.}
\end{minipage} \hskip1cm
\begin{minipage}[t]{8cm}
\centering \includegraphics[width=8cm]{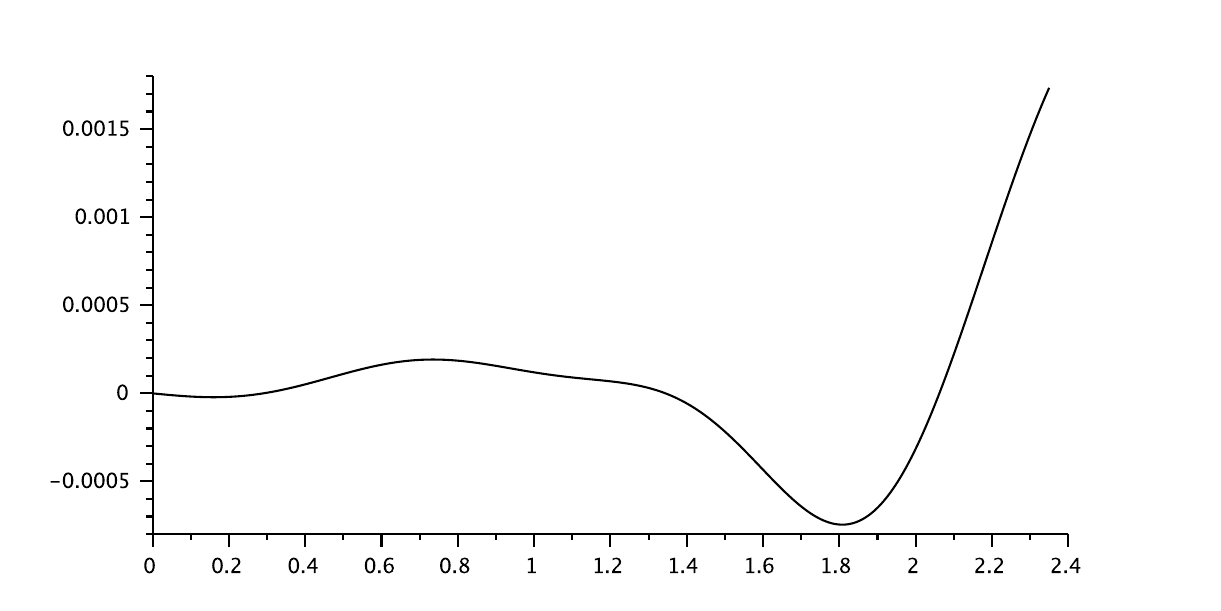} 
\caption{Curve $\frac {d}{d\theta }p(f_1(e^{i\theta }))$, $0\leq \theta \leq \frac{3\pi }{4}$.}
\end{minipage} \end{figure}

b) We turn now to the part $\frac{3\pi }{4}\leq \theta \leq \pi $. We have computed 
 $\displaystyle\min_{750\leq j\leq 1000}\Re(f_1(e^{i\theta_j }))=-0.4998968$ and 
$\displaystyle\max_{750\leq j\leq 1000}\Big|\frac{\Re(f_1(e^{i\theta_{j+1 }}))-\Re(f_1(e^{i\theta_{j }}))}{\theta _{j+1}-\theta _j}\big|=0.01485\dots$ This gives us, for $3\pi /4\leq \theta \leq \pi $, an estimate of $\max\big|\frac {d}{d\theta }p(f_1(e^{i\theta }))\big|\leq 0.015$ and thus $\min \Re(f_1(e^{i\theta}))\geq -0.4998968-0.015\pi /2000>-0.499921$. Note also that this part of the curve clearly satisfies $\Re(f_1(e^{i\theta}))\leq 0$ and $\Im f_1(e^{i\theta}))\leq \sqrt3/6$.
This shows that the set $\{z=f_1(e^{i\theta})\,: \frac{3\pi }{4}\leq \theta \leq \pi \}$ is interior to $W(A)$.
\begin{figure}[!ht] \centering
\begin{minipage}[t]{6cm} \centering
\includegraphics[width=8cm]{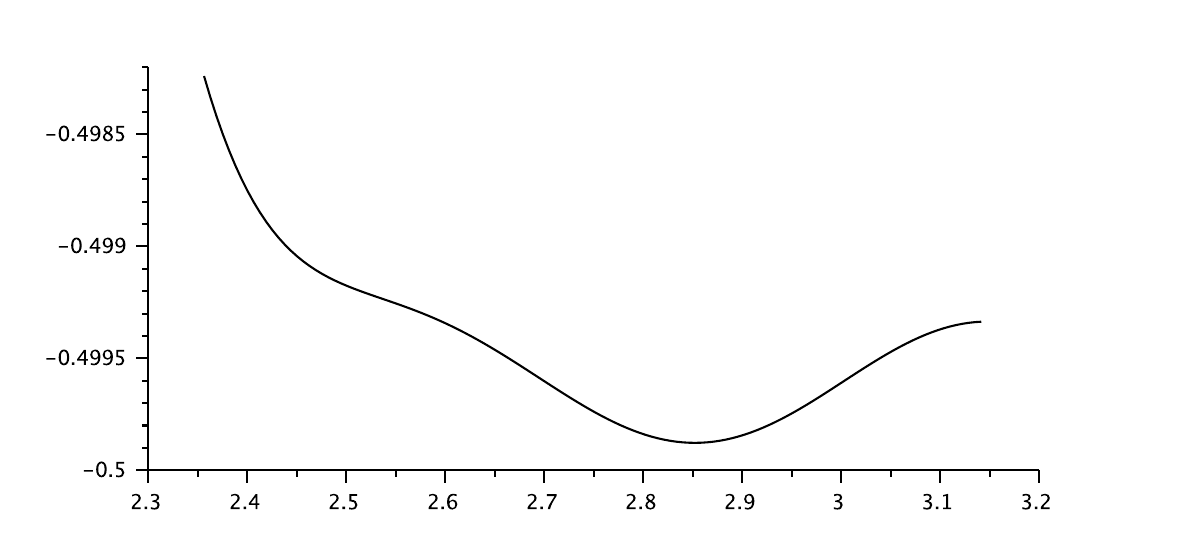}
\caption{\!Curve $\Re(f_1(e^{i\theta}))$,$\frac{3\pi }{4}\leq \theta\! \leq \pi $.}
\end{minipage} \hskip1cm
\begin{minipage}[t]{8cm}
\centering \includegraphics[width=8cm]{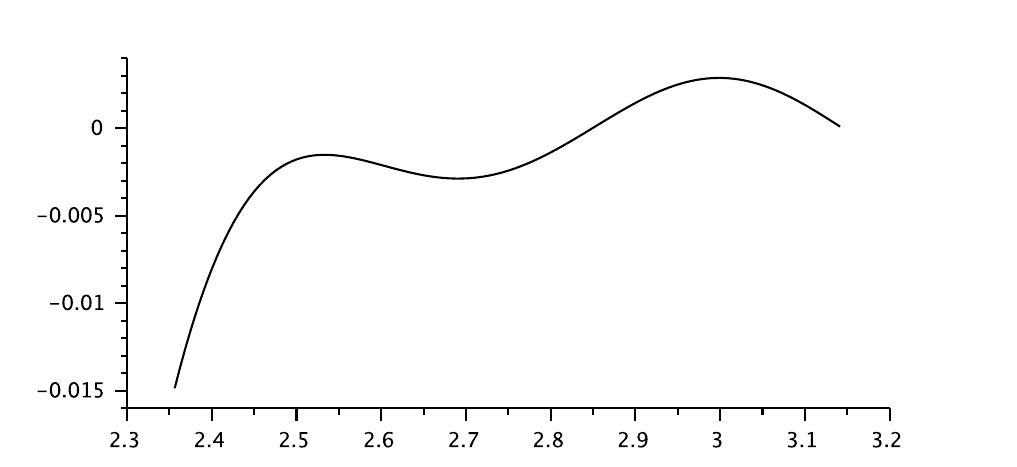} 
\caption{Curve $\frac {d}{d\theta }\Re(f_1(e^{i\theta }))$, $\frac{3\pi }{4}\leq \theta \leq \pi $.}

\end{minipage} \end{figure}

\section{Estimates for the class of matrices $A_k$}

Recall\,\cite{kip2} that the boundary of the  numerical range $W(A)$ is the convex hull of the algebraic curve with tangential equation $T(u,v,w):=$ det$(uB{+}vC{+}wI)=0$ where we have written $A=B{+}i\,C$, with $B$ and $C$ self-adjoint. For the matrix $A_k$, the corresponding tangential equation $T_k(u,v,w)=0$, can be obtained from the recursion
\[
T_1(u,v,w)=w,\quad T_{k+1}(u,v,w)=w\,T_k(u,v,w)+\frac{(w-\frac{u+iv}{2})^k-(w-\frac{u-iv}{2})^k}{i\,v}\frac{u^2+v^2}{4}.
\]
For instance,
\begin{align*}
T_3(u,v,w)&=w^3-\frac34 w(u^2{+}v^2)+\frac14(u^2{+}v^2)u,\\
T_4(u,v,w)&=w^4-\frac32 w^2(u^2{+}v^2)+w(u^2{+}v^2)u-\frac{1}{16}(u^2{+}v^2)(3u^2{-}v^2).
\end{align*}
We can see, by recursion, that
\begin{align*}
T_k(\cos\varphi ,\sin\varphi ,w)=\frac{(-1)^k}{\sin\varphi }\Im\big(e^{-i\varphi }(\tfrac12 e^{i\varphi}{-}w )^k\big).
\end{align*}
We now remark that, if  $\varphi =\frac{k\theta }{2}$ and $w=-\frac12\frac{\sin((k{-}1)\theta /2)}{\sin{(\theta /2})}$, then $e^{-i\varphi /k}(\tfrac12e^{i\varphi }-w)=\frac12\frac{\sin(k\theta /2)}{\sin\theta /2}\in\R$, whence
$T_k(\cos\varphi ,\sin\varphi ,w)=0$.
Now, we consider the algebraic curve $\{ f_k(e^{i\theta })\,: |\theta | \leq \pi \}$ with $f_k(z)=(\sum_{j=1}^{k-1}jz^{k-j})/k$ $=\frac{z^{k+1}-z-k(z^2-z)}{k(z-1)^2}$. We remark that
\[
k\,e^{i\theta}\,f'_k(e^{i\theta })=\sum_{j=1}^{k-1}j(k{-}j)e^{i(k-j)\theta} =\frac{e^{ik\theta /2}}{2}\sum_{j=1}^{k-1}j(k{-}j)\cos\tfrac{(k{-}2j)\theta}{2};
\]
hence, the vector $e^{ik\theta /2}$ is the unit normal at the point $ f_k(e^{i\theta }) $ and the equation of the tangent at this point is $\Re\big((e^{-ik\theta/2}(x{+}iy{-}f(e^{i\theta })\big) =0$, i.e.
$ux{+}vy{+}w=0$, with $u=\cos\frac{k\theta }2$, $v=\sin\frac{k\theta }2$, 
$w=-\Re \big(e^{-ik\theta /2}f(e^{i\theta })\big)=-\frac12\sum_{j=1}^{k-1}\cos(\frac{k-2j}2\theta)=-\frac12\frac{\sin((k-1)\theta /2)}{\sin(\theta /2)} $. Thus, this shows that
 the algebraic curve which satisfies the tangential equation $T_k(u,v,w)=0$ is the set:\footnote{Thanks to Bernd Beckermann for a useful remark.}

$\{z\,: z=\frac1k\sum_{j=1}^{k-1}\,j\,e^{i(k-j)\theta },\quad-\pi \leq \theta \leq \pi \}$.

\begin{figure}[!ht] \centering
\begin{minipage}[t]{4cm} \centering
\includegraphics[width=4cm]{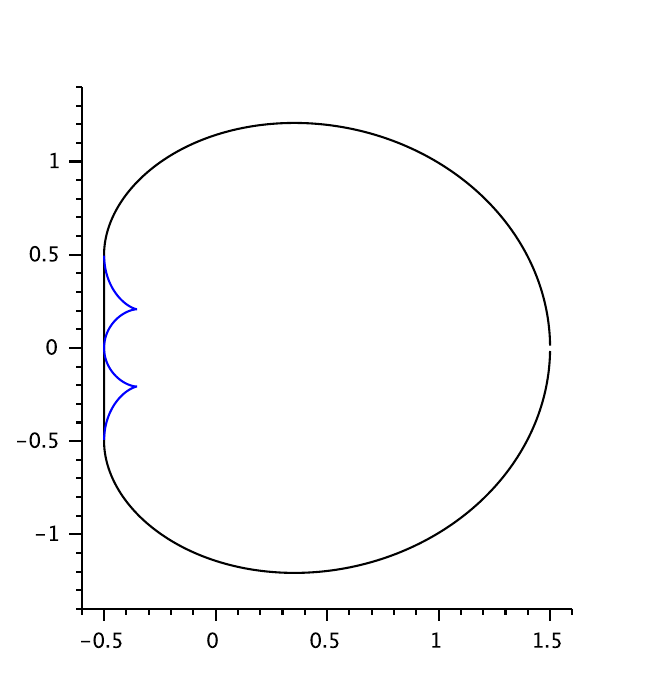}
\caption{$W(A_4)$.}
\end{minipage} \hskip1cm
\begin{minipage}[t]{4cm} \centering
\includegraphics[width=4cm]{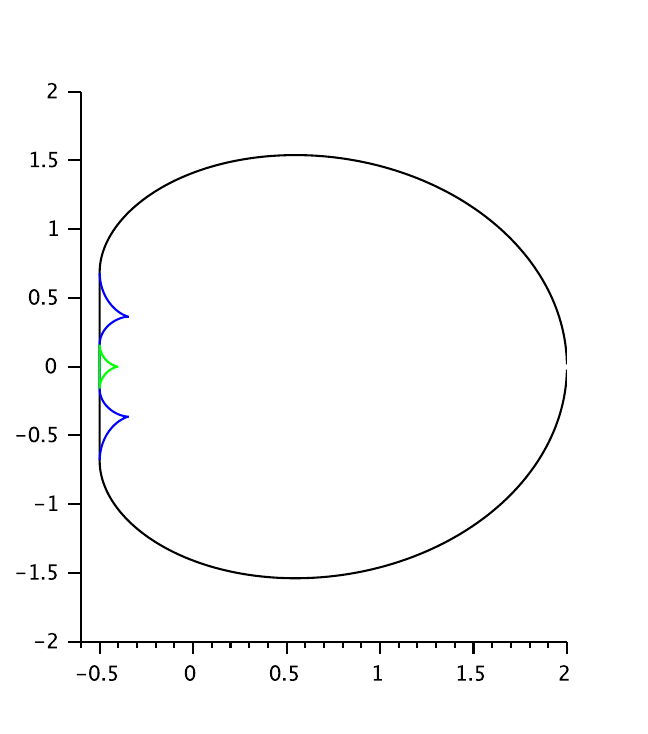}
\caption{$W(A_5)$.}
\end{minipage} \hskip1cm
\begin{minipage}[t]{4cm}
\centering \includegraphics[width=4cm]{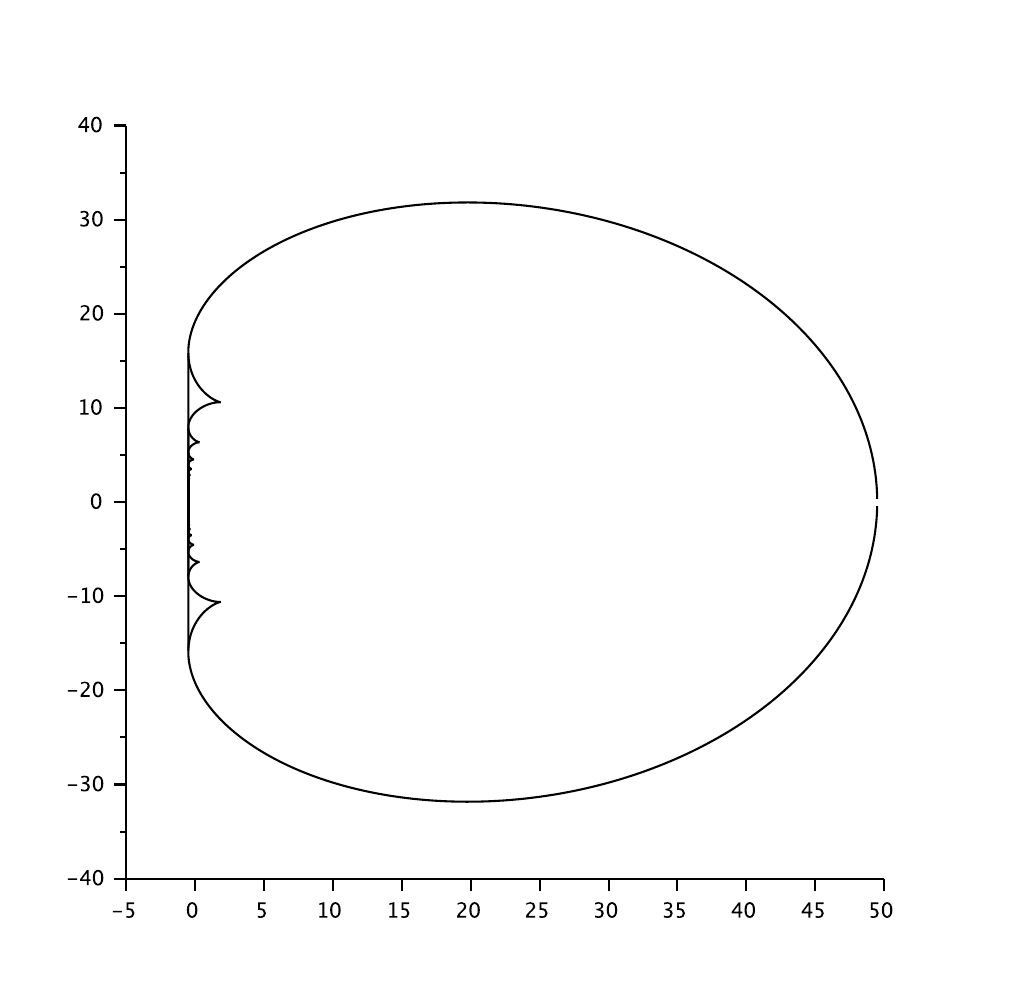} 
\caption{$W(A_{100})$.}
\end{minipage} \end{figure}

The points with horizontal tangent are given by $\theta _j=\frac{(2j-1)\pi }{k}$, $j=1,\dots,n$ and are cuspid for $j=2,\dots,k{-}1$; the $k{-}1$ points of the algebraic curve on the flat part 	are the points $-\frac12{+}\frac{i}{2\tan(j \pi /k)}$, $j=1,\dots,k{-}1$. \medskip

Note that the matrix $A_k{+}A_k^*{+}I$, whose entries are all $1$'s, has the simple eigenvalue $k$ (eigenvector $(1,1,\cdots,1)^T$),
and the eigenvalue $0$ of multiplicity $k{-}1$ (with orthonormal eigenvectors $\frac{1}{\sqrt{k}}(e^{i\theta _j},e^{2i\theta _j},\dots,e^{ik\theta _j})^T$, $\theta _j=\frac{2j\pi }{k}$, $j=1\dots,k{-}1$). This implies that the point $\frac{k-1}{2}=f_k(1)$ is the extremal right point of $W(A_k)$, and that
the boundary has a flat part on the line $\Re z=-\frac12$.  The boundary of the numerical range is the union of a part of the algebraic curve and of a straight part:
\[
\partial W(A_k)=\{z: z=\tfrac1k\sum_{j=1}^{k-1}j\,e^{i(k-j)\theta }, |\theta| \leq \tfrac{2\pi}k \}\cup
\{z: z=-\tfrac12(1{+}iy),\ -\cot(\tfrac\pi k) \leq y\leq \cot(\tfrac\pi k) \}.
\] 

\bigskip
We turn now to some estimates for the constant corresponding to the matrices $A_k$, $k\leq 6$. 
For $A_4$, we have computed the values $g'(0)=1.1888506$, $g''(0)=-1.6292742$, $g'''(0)=4.7085601$, which gives $g(A_4)=\begin{pmatrix}0 &a &b &c\\0&0 &a &b\\0 &0 &0&a\\0&0&0&0\end{pmatrix}$, with
$a=1.1888506$, $b=0.3742134$, $c=0.3443362$.
 
We have obtained a lower bound $\psi (A_4)\geq 1.9938003$ with the Blaschke product corresponding to the 3 coefficients $-0.4560323 {\pm} 0.3891911i$, $0.2474013$. 
We have also an upper bound: with the values $x=b/a^{1.5}$, $y=az{-}bx/a{+}c/a^{1.5}$, $z=-0.0735033$,
$t=-0.0231366$ and the matrix\\[8pt]
$H=\begin{pmatrix}a^{1.5 }&xa &y &t\\0&a^{0.5} &0 &z\\0 &0 &a^{-0.5}&-xa\\ 0&0&0&a^{-1.5}\end{pmatrix}$, it holds  $H^{-1}g(A_4)H=\begin{pmatrix}0 &1 &0 &0\\0&0 &1 &0\\0 &0 &0&1\\0&0&0&0\end{pmatrix}$.\\ Then, cond$(H)\simeq1.9938002$ and $\|H^{-1}g(A_4)H\|=1$. 
 
We can consider that the estimate $1.993800\leq \psi (A_4)\leq 1.993801$ is correct.\bigskip
 
For $A_5$, we have computed  $g(A_5)=\begin{pmatrix}0 &a &b &c &d\\0&0 &a &b &c\\0 &0 &0&a& b\\0&0&0&0 &a\\0 &0 &0 &0 &0\end{pmatrix}$,
with $a=1.1170233$, $b=0.2325756$, \\$c=0.2187502$, $d=0.1895824$. 

We have obtained a lower bound $\psi (A_5)\geq 1.9929216$ with the Blaschke product corresponding to the 4 coefficients $-0.2583004 {\pm} 0.60451151i$, $-0.6247827$, $0.3295365$. 
We have also an upper bound: with the values $u=-0.0194597$, $w=-0.0384976$, $g=-0.1091772$,
$h=-0.2503045$, $f=ah+ba^{-2}$, $v=ag{+}bh+ca^{-2}$, $z=af{+}ba^{-1}$, $t=aw{+}bg{+}ch{+}da^{-2}$
, $y=av{+}bf{+}ca^{-1}$, $x=az{+}b$ and the matrix\\[8pt]
 $H=\begin{pmatrix}a^2&x &y &t &u\\0&a &z &v&w\\0 &0 &1&f&g\\ 0&0&0&a^{-1} &h\\ 0&0&0&0&a^{-2}\end{pmatrix}$, it holds  $H^{-1}g(A_5)H=\begin{pmatrix}0 &1 &0 &0&0\\0&0 &1 &0&0\\0 &0 &0&1&0\\0&0&0&0&1\\0 &0 &0 &0 &0 \end{pmatrix}$. \\Then, cond$(H)\simeq1.9929216$ and $\|H^{-1}g(A_5)H\|=1$. 
 
 We can consider that the estimate $1.992921\leq \psi (A_5)\leq 1.992922$ is correct.\bigskip

For $A_6$, we have computed  $g(A_6)=\begin{pmatrix}0 &a &b &c &d &e\\0&0 &a &b &c &d\\0 &0 &0&a& b&c\\0&0&0&0 &a &b\\0 &0 &0 &0 &0 &a\\0 &0 &0 &0 &0 &0 \end{pmatrix}$,
with $a=1.0798634$, $b=0.1590093$, $c=0.1519169$, $d=0.1359021$ et $e=0.1161184$. 

We have obtained a lower bound $\psi (A_6)\geq 1.9924447$ with the Blaschke product corresponding to the 5 coefficients $-0.5859775 {\pm} 0.3199164 i$,  $-0.0604565 {\pm} 0.70030221 i$,  $0.3972632$, $0.3295365$. 
We have also a upper bound: with the values \\
$y=(-0.0163999,-0.0248879, -0.0578414,-0.1294105,-0.243031)$,\\
$x_9=a\,y_5+b\,a^{-2.5}$,  $x_8=a\,y_4+b\,y_5+c\,a^{-2.5}$,
$x_7=a\,x_9+b\,a^{-1.5}$, $x_6=a\,y_3+b\,y_4+c\,y_5+d\,a^{-2.5}$,
$x_5=a\,x_8+b\,x_9+c\,a^{-1.5}$, $x_4=a\,x_7+b\,a^{-0.5}$, $x_3=a\,y_2+b\,y_3+c\,y_4+d\,y_5+e\,a^{-2.5}$,\\
$x_2=a\,x_6+b\,x_8+c\,x_9+d\,a^{-1.5}$, $x_1=a\,x_5+b\,x_7+c\,a^{-0.5}$; $x_0=a\,x_4+b\,a^{0.5}$, 

 \noindent and the matrix\\[5pt]
 $H=\begin{pmatrix}a^{2.5}&x_0 &x_1&x_2&x_3&y_1\\0&a^{1.5} &x_4 &x_5&x_6&y_2\\0 &0 &a^{.5}&x_7&x_8&y_3\\ 0&0&0&a^{-.5} &x_9 &y_4\\ 0&0&0&0&a^{-1.5}&y_5\\0 &0 &0 &0 &0 &a^{-2.5}\end{pmatrix}$, it holds  $H^{-1}g(A_6)H=\begin{pmatrix}0 &1 &0 &0&0&0\\0&0 &1 &0&0&0\\0 &0 &0&1&0&0\\0&0&0&0&1&0\\0 &0 &0 &0 &0&1\\0 &0 &0 &0 &0 &0 \end{pmatrix}$. \\Then, cond$(H)\simeq1.9924445$ and $\|H^{-1}g(A_6)H\|=1$. 
 
 We can consider that the estimate $1.992444\leq \psi(A_6)\leq 1.992445$ is correct.

\noindent{\bf Remark}. {\it In each of these cases, the matrix $M=g(A_k)$ satisfies the relation
$\psi _{\D}(M)=\psi _{cb,\D}(M)$. This property holds for all  $d\times d$ matrices $M$ if $d\leq 2$, but may fail {\rm\cite{cgh}} if $d\geq 3$. Also here, for the matrix $H$ which realizes $\psi _{cb,\D}(M)$, 
$H^{-1}MH$ was a Jordan block, which is not generally the case}.

Table \ref{table1} summarizes our results.
\begin{table}[h]
\begin{center}
\begin{tabular}{|c|c|c|c|} \hline
$n$ & lower bound & upper bound & difference \\ \hline
$3$ & $1.9956978$ & $1.9956979$ & $10^{-7}$ \\ \hline
$4$ & $1.993800$  & $1.993801$  & $10^{-6}$ \\ \hline
$5$ & $1.992921$  & $1.992922$  & $10^{-6}$ \\ \hline
$6$ & $1.992444$  & $1.992445$  & $10^{-6}$ \\ \hline
\end{tabular}
\end{center}
\caption{Upper and Lower Bounds on $\psi ( A_n ) = \psi_{cb,\D}( M_n )$. \label{table1}}
\end{table}

\section{ About the computation of the conformal mapping $g$ }

We may write $g(z)=z\,\exp(u\!+\!iv)$, with $u(z)$ and  $v(z)$ harmonic real-valued functions.
Note that $u(z)=-\log |z| $ on $\partial W(A)$, which determines $u$ in $W(A)$ in a unique way.

\smallskip

Let us consider a representation
$\partial W(A)=\{\sigma(\theta)\,; \theta\in [0,2\pi]\}$ of the boundary and choose $\lambda>0$.
If $q$ is a $2\pi$-periodic real-valued function such that,  
\begin{align*}
\int_{0}^{2\pi} q(\theta)\log\big|\frac{\sigma(\theta){-}\sigma(\varphi)}{\lambda }\big|\, d\theta= -\log|\sigma(\varphi)|,\quad
\text{ for all }\varphi\in[0,2\pi[,
\end{align*}
then it holds $u(z)=\int_{0}^{2\pi} q(\theta)\log\big|\frac{\sigma(\theta){-}z}{\lambda }\big|\, d\theta$
since this integral is clearly harmonic and is equal to $-\log |z| $ on $\partial W(A)$. It is known that
such a $q$ exists if and only if $\lambda $ is different from the logarithmic capacity of $W(A)$.
Generally we will use this equation with $\lambda =1$ and then rewrite it as
\[
\int_{0}^{2\pi} q(\theta)\log\Big|\frac{\sigma(\theta){-}\sigma(\varphi)}{e^{i\theta}{-}e^{i\varphi}}\Big|\, d\theta +
\int_{0}^{2\pi} \!\!q(\theta)\log |e^{i\theta}{-}e^{i\varphi}|\,  d\theta= -\log|\sigma(\varphi)|,\ \ 
\forall\,\varphi\in[0,2\pi[.
\]
We discretized this equation using a representation  $\sigma(\theta)$ of $\partial W(A)$ and approximating $q(\cdot)$ by a trigonometric polynomial $q_{n}(\cdot)$ of degree $n$, and employing a collocation method at the points $\theta_{j}$, $j=0,1,\dots, 2n$ (it is known that an odd number of collocation points is necessary  for such a method). So, we get an approximation $q_{j}=q_{n}(\theta_{j})$ by solving the system
\begin{align*}
\frac{2\pi}{2n+1}
\sum_{j=0}^{2n} q_{j} \log\Big|\frac{\sigma(\theta_{j})-\sigma(\theta_{i})}{e^{i\theta_{j}}-e^{i\theta_{i}}}\Big|+
\int_{0}^{2\pi} q_{n}(\theta)\log |e^{i\theta}-e^{i\theta_{i}}|\, d\theta= -\log|\sigma(\theta_{i})|,\\
\text{ for }i=0,1,\dots, 2n.
\end{align*}
We have approximated the first integral by the trapezoidal formula; of course, if $j=i$, we have to replace $\log\Big|\frac{\sigma(\theta_{j})-\sigma(\theta_{i})}{e^{i\theta_{j}}-e^{i\theta_{i}}}\Big|$ by 
$\log |\sigma'(\theta_{i})|$.
 Recall that, for the remaining integral, there holds
\begin{align*}
\int_{0}^{2\pi} q_{n}(\theta)\log |e^{i\theta}-e^{i\theta_{i}}|\, d\theta= -\frac{2\pi}{2n+1}\sum_{j=0}^{2n} c(j\!-\!i) \,q_{j},\\
\text{with  }c(k)=c(-k)=\sum_{j=1}^{n}\frac{\cos j\theta_{k}}{j}.
\end{align*}
Then, we obtain the approximation of $u$ from
\[
u(z)\simeq\frac{2\pi}{2n+1}
\sum_{j=0}^{2n} q_{j} \log|\sigma(\theta_{j})-z|
\]
and the approximation of the derivatives of $g$ at $0$ (note that here $v(x)=0$ if $x\in \R$)
\begin{align*}
g'(0) &= \exp(u(0)), \ g''(0)=2\,g'(0)u'(0),\ g^{(3)}(0)=3\,g'(0)(u'(0)^2{+}u''(0)),\\
g^{(4)}(0) &= 4\,\,g'(0)(u'(0)^3{+}3u'(0)u''(0){+}u^{(3)}(0)),\\
g^{(5)}(0)&= 5\,\,g'(0)(u'(0)^4{+}6u'(0)^2u''(0){+}4u'(0)u^{(3)}(0){+}3u''(0)^2{+}u^{(4)}(0)),
\end{align*}
via the formulae
\[
g'(0)\simeq\exp\Big(\frac{2\pi }{2n+1}\sum_{j=1}^{2n+1}q_j\log|\sigma (\theta _j)|\Big),
\quad  u^{(k)}(0)\simeq-(k{-}1)!\frac{2\pi }{2n+1}\Re\Big(\sum_{j=1}^{2n+1}\frac{q_j}{\sigma (\theta _j)^k})\Big).
\]

\noindent{\bf Remark.} \textit{ In order to get $q$, we have to solve a linear system of the form $Mq=b$. But it could appear that the matrix $M$ is not invertible, or badly conditioned, if the logarithmic capacity of $W(A)$ is close to 1. In this case, we can replace this system by $(M{-}E)q=b{-}e$
where $E$ (resp. $e$) is a matrix (resp. a vector) with all entries equal to $1$.}

This method is very efficient for analytic boundary (exponential convergence with respect to $n$, see for instance \cite{ChAr}), but here we have singularities at the transition points between the straight line and the algebraic part, which reduces the order of convergence to $O(n^{-4})$ which is still good. With
$\varphi_k(z)=\tfrac1k\sum_{j=1}^{k-1}j\,z^{i(k-j)}$, we have
\[
\partial W(A_k)=\{z: z=\varphi _k(\theta), -\tfrac{2\pi}k\leq \theta \leq \tfrac{2\pi}k \}\cup
\{z: z=-\tfrac12(1{+}iy),\ -\cot(\tfrac\pi k) \leq y\leq \cot(\tfrac\pi k) \}.
\] 
We have used $2n{+}1$ points on the algebraic part of $\partial W(A_k)$: $z_j= \varphi _k(\frac{2\pi j}{kn})$ for $j=-n,\dots,n$ and $2n_2$ equidistant points on the straight part $z_{n+j}=z_n{-}jh\,i$, $j=1,\dots,2n_2$ where $h=2\cot(\tfrac\pi k)/(2n_2{+}1)$ and $n_2$ is chosen such that $h$ is as close as possible to $|z_n{-}z_{n-1}|$.

\section{ Program in Matlab for the computation of $g(A)$} 

\noindent\texttt{function [gofA,gderivs,nn]= Akstudy(k,n);\\
\\
\textcolor{blue}{\%  For 3 <= k <= 6, forms the kxk matrix A with ones in the\\
 \% strict upper triangle and zeros elsewhere, and computes\\
\% its image gofA under the Riemann mapping from W(A) to the\\
\%  unit disk with g(0) = 0, g'(0) > 0.  gderivs(j), j=1,...,5,\\
\%  contains the value of the jth derivative of g at 0.\\
\% W(A) consists of a cardioid and a vertical line segment,\\
\% and 2n+1 points are used to represent the cardioid portion.\\
\% Output argument nn is then the total number of discretization\\
\%  points used to represent the boundary of W(A).  Note also\\
\% that n may be modified (to make things come out even).\\
\\
}A = triu(ones(k),1);  \textcolor{blue}{\% Form the matrix.}}\\

\noindent\texttt{\textcolor{blue}{\%  Discretize W(A).  Use n+1 points on the upper part of the\\ 
\%  algebraic curve}}.

\noindent\texttt{th = 2*[0:n]'*pi/n/k;\\
z = zeros(size(th));\\
for j=1:k-1, z = z + (k-j)*exp(1i*j*th)/k; end;}\\

\noindent\texttt{\textcolor{blue}{\%   Choose the same step size on the line segment.}\\
h = abs(z(n+1)-z(n));\\
nn = fix(imag(z(n+1))/h - 0.5);\\
h = imag(z(n+1))/(nn+0.5);\\
for j=1:nn, z(n+1+j) = z(n+1) - 1i*j*h; end;}\\

\noindent\texttt{\textcolor{blue}{\%    Complete by symmetry.}\\
nn = n+1+nn;\\
zz = conj(z);\\
z = [z(1:nn); zz(nn:-1:2)];}\\

\noindent\texttt{\textcolor{blue}{\%    Plot W(A).}\\
plot([z; z(1)],'-k','LineWidth',2), axis equal, shg}\\

\noindent\texttt{\textcolor{blue}{\%  Compute the conformal mapping.}\\
nn = length(z); n = (nn-1)/2;\\
zz = z;\\
e = [1:n]'; ee = 2*pi/nn * [1:n]';\\
c0 = sum(ones(size(e)) ./ e);\\
c = zeros(nn-1,1); d = zeros(nn-1,1);\\
for j=1:nn-1, c(j) = sum(cos(j*ee)./e); end;\\
for j=1:nn-1, d(j) = sum(2*sin(j*ee) .* e)/nn; end;\\
dd = [d; 0];\\
zzprim = zeros(nn,1);\\
for j=1:nn, dd = [dd(nn); dd(1:nn-1)]; zzprim(j) = sum(zz.*dd); end;\\
zzprim = abs(zzprim);}\\

\noindent\texttt{\textcolor{blue}{\% Compute the matrix M such that Mq = -log |sigma|.}\\
ee = exp(1i*2*pi/nn * [1:nn]');\\
M = zeros(nn,nn);\\
for j=1:nn,\\
\phantom{1}\   M(j,j) = log(zzprim(j)) - c0;\\
\phantom{1}\   for k=j+1:nn,\\
\phantom{1}\ \  \    M(k,j) = log(abs((zz(k)-zz(j))/(ee(k)-ee(j)))) - c(k-j);\\
\phantom{1}\ \  \    M(j,k) = M(k,j);\\
\phantom{1}\  end;\\
end;}\\

\noindent\texttt{\textcolor{blue}{\%  If t < 0, M is badly conditioned, so we translate M.}\\
t = 10$\scriptsize\wedge$4-cond(M);\\
if t < 0, M = M - ones(M); t, pause, end;\\}

\noindent\texttt{\textcolor{blue}{\%  Compute q.}\\
b0 = log(abs(zz));\\
q = -M$\backslash$b0;}\\

\noindent\texttt{\textcolor{blue}{\%  Take account of the translation.}\\
if t < 0, b0 = b0 - ones(b0); end;}\\

\noindent\texttt{\textcolor{blue}{\%  Compute derivatives of g at 0.}\\
gp = exp(sum(q.*b0));\\
b = -real(sum(q./zz));\\
c = -real(sum(q./zz.$\scriptsize\wedge$2));\\
d = -2*real(sum(q./zz.$\scriptsize\wedge$3));\\
ed = -6*real(sum(q./zz.$\scriptsize\wedge$4));\\
gs = 2*gp*b; gt = 3*gp*(b$\scriptsize\wedge$2 + c);\\
gq = 4*gp*(3*b*c + b$\scriptsize\wedge$3 + d);\\
gc = 5*gp*(b$\scriptsize\wedge$4 + 6*b$\scriptsize\wedge$2*c + 4*b*d + 3*c$\scriptsize\wedge$2 + ed);\\
gderivs = [gp, gs, gt, gq, gc];\\
gofA = gp*A + 0.5*gs*A$\scriptsize\wedge$2 + (1/6)*gt*A$\scriptsize\wedge$3 + (1/24)*gq*A$\scriptsize\wedge$4 + (1/120)*gc*A$\scriptsize\wedge$5;}\\

\bigskip
\noindent {\bf Acknowledgment:} The first author conducted this work within the France 2030 framework programme, the Centre Henri Lebesgue  ANR-11-LABX-0020-01


\begin{thebibliography}{11}

\bibitem{ChAr}  {\sc R.S.\,Cheng and D.\,Arnold}, {\em The delta-trigonometric method using the single layer potential representation},
J. of Integral Equations and Applications, {\bf 1} (1988), pp.~517--547, \verb+https://doi.org/10.1216/JIE-1988-1-4-517+.

\bibitem{Choi2013}{\sc D.~Choi}, {\em A proof of Crouzeix's conjecture for 
a class of matrices}, {Lin.~Alg.~Appl.}, 438 (2013), pp.~3247--3257, \verb+https://doi.org/10.1016/j.laa.2012.12.045+.

\bibitem{ChoiGreen2015}{\sc D.~Choi, A.~Greenbaum}, {\em Roots of matrices in the study of GMRES
convergence and Crouzeix's conjecture}, {SIAM J.~Matrix~Anal.~Appl.}, 36 (2015), pp.~289--301,
\verb+https://doi.org/10.1137/140961742+.

\bibitem{Crouzeix2004}{\sc M.~Crouzeix}, {\em Bounds for analytic functions 
of matrices}, {Int.~Eq.~Op.~Th.}, 48 (2004), pp.~461-477,
\verb+https://doi.org/10.1007/s00020-002-1188-6+.


\bibitem{cgh} {\sc M.\,Crouzeix, F.\,Gilfeather, and J.\,Holbrook}, {\em Polynomial bounds for small matrices},
    {Linear and Multilinear Algebra}, {\bf 62}, no.~5 (2014), pp.~614--625, \verb+https://doi.org/10.1080/03081087.2013.777439+.
    
\bibitem{GKL2018}{\sc C.~Glader, M.~Kurula, and M.~Lindstrom},
{\em Crouzeix's conjecture holds for tridiagonal $3 \times 3$ matrices
with elliptic numerical range centered at an eigenvalue}, 
{SIAM J.~Matrix~Anal.~Appl.}, 39 (2018), pp.~346--364, \verb+https://doi.org/10.1137/17M1110663+.

\bibitem{GW2013}{\sc H.-L.~Gau and P. Y.~Wu}, {\em Numerical ranges of KMS matrices}, Acta Sci. Math. (Szeged), 79 (2013), no. 3--4, pp.~583--610,
\verb+https://doi.org/10.1007/BF03651342+.

\bibitem{KMS1953}{\sc M.~Kac, W. L.~Murdock and G.~Szeg\H o}, {\em On the eigenvalues of certain Hermitian forms}, J. Rational Mech. Anal., 2 (1953), pp.~767--800,
\verb+https://www.jstor.org/stable/24900353+.

\bibitem{kip2}{\sc R.\,Kippenhahn}, {\em On the numerical range of a matrix}, Translated from the German by Paul F. Zachlin and Michiel E. Hochstenbach, {Linear Multilinear Algebra} {\bf 56} (2008), pp.~185--225, \verb+https://doi.org/10.1080/03081080701553768+.


\bibitem{Overton2022} {\sc M.~L.~Overton}, {\em Local minimizers of the Crouzeix
ratio:  a nonsmooth optimization case study}, {Calcolo} {\bf 59,8} (2022),
\verb+https://doi.org/10.1007/s10092-021-00448-z+.

  
 \bibitem{paul} {\sc V.\,Paulsen},
{\em Completely bounded maps and operator algebras}, Cambridge Univ. Press, 2002,
\verb+ https://doi.org/10.1017/CBO9780511546631+.



\end{thebibliography}
\end{document}